\documentclass[11pt]{paper}

\usepackage{amssymb,amsmath,enumerate,theorem,epsfig,rotating,graphics,changebar,eepic}
\usepackage{amsfonts}
\usepackage{a4,latexsym,parskip}
\usepackage{hyperref}
\hypersetup{pdfauthor=MS} \hypersetup{pdftitle=Fields}

\newcommand{\PP}{\mathbb{P}}

\newcommand{\C}{\mathbb{C}}
\newcommand{\Q} {\mathbb{Q}}

\newcommand{\F}{\mathbb{F}}
\newcommand{\Z}{\mathbb{Z}}

\newcommand{\OO}{\mathcal{O}}

\newcommand{\jj}{\mathfrak{j}}

\newcommand{\qed}{\hfill \ensuremath{\Box}}
\newcommand{\NS}{\mathop{\rm NS}\nolimits}
\newcommand{\A}{\mathfrak{A}}
\newcommand{\TT}{\mathfrak{T}}
\newcommand{\Aut}{\mathop{\rm Aut}\nolimits}

\theoremstyle{break} \newtheorem{Theorem}{Theorem}
\newtheorem{Proposition}[Theorem]{Proposition}
\newtheorem{Lemma}[Theorem]{Lemma}
\newtheorem{Example}[Theorem]{Example}

\newtheorem{Definition}[Theorem]{Definition}
\newtheorem{Corollary}[Theorem]{Corollary}

\newtheorem{Remark}[Theorem]{Remark}

\begin{document}
\setlength{\unitlength}{1cm}

\title{K3 surfaces with non-symplectic automorphisms of 2-power order} 

\author{Matthias Sch\"utt}


\date{\today}
\maketitle

\abstract{This paper concerns complex algebraic K3 surfaces with an automorphism which acts trivially on the N\'eron-Severi group. Complementing a result by Vorontsov and Kond\=o, we determine those K3 surfaces where the order of the automorphism is a 2-power and equals the rank of the transcendental lattice. We also study the arithmetic of these K3 surfaces and comment on relations with mirror symmetry.}




\section{Introduction}
\label{s:intro}

This paper concerns complex algebraic K3 surfaces with a non-symplectic automorphism which acts trivially on the algebraic cycles. 
Such K3 surfaces have been extensively studied using lattice theory introduced by Nikulin. The classification of these K3 surfaces due to Vorontsov and Kond\=o is twofold (cf.~
Theorems~\ref{Thm:uni}, \ref{Thm:non}). First it gives all possible orders of the non-symplectic automorphism in general. Then it determines unique K3 surfaces in the extreme case where the transcendental lattice $T(X)$ is as small as possible relative to the order of the automorphism -- but only for orders which are not powers of 2. This paper complements the results of Vorontsov and Kond\=o by virtue of the following classification:

\begin{Theorem}\label{thm}
Let $X$ be a K3 surface with a non-symplectic automorphism $\eta$ which acts trivially on $\NS(X)$. Assume that the order $m$ of $\eta$ is a $2$-power and that $T(X)$ has  rank $m$. Up to isomorphism we are in one of the following cases:
$$
\begin{array}{c|c|c|c}
\hline
m & \NS(X) & \text{equation} & T(X)\\
\hline \hline
2 & U+E_8^2+A_1^2 & y^2 = x^3 - 3\,t^4\,x + t^5 + t^7 & \langle 2\rangle^2\\
\hline
4 & U+E_8^2 & y^2 = x^3 - 3\,\lambda\,t^4\,x+t^5+t^7 & U^2\\
\hline
4 & U+D_8+E_8 & y^2 = x^3 + t\,x^2 + \lambda\,t^4 x+t^7 & U+U(2)\\
\hline
8 & U+D_4+E_8 & y^2 = x^3 +\lambda\,t\,x^2+t^2\,x+t^7 & U^2+D_4\\
\hline
16 & U+D_4 & y^2 = x^3+\lambda\,t\,x^2+t^2\,x+t^{11} & U^2+D_4+E_8\\
\hline
\end{array}
$$
In each case, a general choice of the parameter $\lambda$ guarantees that $T(X)$ really does have rank $m$. 
\end{Theorem}

Here $U$ denotes the hyperbolic plane with intersection form $\begin{pmatrix}0 & 1\\1 & 0\end{pmatrix}$, and $U(2)$ indicates that the intersection form is multiplied by $2$.
The lattices $A_n, D_k, E_l$ refer to negative-definite simple root lattices;
these are in correspondence with Dynkin diagrams.
The lattice $\langle 2\rangle$ is generated by a single element of self-intersection $2$.
Thus $\langle 2\rangle$ and $A_1$ agree up to the sign of the intersection form.

Note that only the first case for $m=4$ is unimodular. In all other cases, $\NS(X)$ has discriminant $-4$. Details concerning the general choice of $\lambda$ can be found in (\ref{eq:OO}) and Corollary~\ref{cor}.

The proof of Theorem~\ref{thm} is based on lattice theory as developed by Nikulin. The main ideas go back to Kond\=o. He used special properties of elliptic fibrations for the classification. We will recall the main concepts in the next section. 
This will culminate in a list of all theoretically possible N\'eron-Severi lattices (Tab.~\ref{tab}). We first consider the unimodular case in section \ref{s:uni}. The main part of this paper is devoted to the non-unimodular cases. Sections \ref{s:non-4}--\ref{s:16} will rule out all lattices but the ones in Theorem~\ref{thm}. We then derive the given families of K3 surfaces.

After the proof of Theorem~\ref{thm}, the paper continues with a discussion of arithmetic aspects. Within the families of Theorem~\ref{thm}, we find K3 surfaces of CM type and determine their zeta functions over finite fields (Theorem~\ref{Thm:zeta}). This result makes use of coverings by Fermat surfaces which we briefly review in section \ref{s:Fermat}. We conclude with comments about mirror symmetry. For the families in Theorem~\ref{thm}, we determine mirror partners with comparable arithmetic properties.

\section{Preliminaries}

In this section, we recall the classification result of Kond\=o and Vorontsov.
Furthermore we review the techniques from lattice theory and basics on elliptic surfaces that will be used to prove Theorem \ref{thm}.
Most of these ideas go back to Nikulin and Kond\=o.

\subsection{The classification of Kond\=o and Vorontsov}
\label{s:class}

Let $X$ be a complex  algebraic K3 surface endowed with an automorphism $\eta$ of order $m$. We call $\eta$ non-symplectic if it acts on the holomorphic $2$-form as multiplication by a primitive $m$-th root of unity $\zeta_m$.

The N\'eron-Severi group $\NS(X)$ of $X$ consists of divisors up to algebraic equivalence. For a K3 surface, we can also consider numerical equivalence instead. Through cup-product, $H^2(X,\Z)$ is endowed with the structure of the unique even unimodular lattice of rank $22$ and signature $(3,19)$:
\[
H^2(X,\Z) = U^3 + E_8^2.
\]
Since $\NS(X)=H^2(X,\Z)\cap H^{1,1}(X)$ by Lefschetz' theorem, it inherits the structure of a lattice. Its rank is called the Picard number $\rho(X)$. By the Hodge index theorem, $\NS(X)$ has signature $(1, \rho(X)-1)$.

The transcendental lattice $T(X)$ is the orthogonal complement of $\NS(X)$:
\[
T(X)=\NS(X)^\bot\subset H^2(X,\Z).
\]
It is known that the representation of $\Aut(X)$ on $\NS(X)+T(X)$ is faithful. I.e.~let $O(\NS(X))$ and $O(T(X))$ denote the respective groups of isometries. Then the induced map 
\[
\Aut(X)\to O(\NS(X)) \times O(T(X))
\]
is injective. It follows that any non-trivial automorphism that acts trivially on $\NS(X)$ is non-symplectic. Another important consequence is that 
\begin{eqnarray}\label{eq:phi}
\phi(m) \,|\, \mbox{rank}(T(X))
\end{eqnarray}
where $m$ is the order of the non-symplectic automorphism $\eta$ and $\phi$ is Euler's $\phi$-function.
This follows from the $\Z[\zeta_m]$-module structure on $T(X)$ given by $\eta$ (cf.~\cite[Theorem~3.1]{Ni}). 
Vorontsov \cite{Vo} announced a classification including all those cases where we have equality in (\ref{eq:phi}). Kond\=o \cite{Ko} corrected and proved the statements.

\begin{Theorem}[unimodular case]\label{Thm:uni}
Let $X$ be an algebraic K3 surface with an automorphism $\eta$ of order $m$. Assume that $\NS(X)$ is unimodular and $\eta$ acts trivially on $\NS(X)$. Let $\Omega=\{12,28,36,42,44,66\}$. Then
\begin{enumerate}[(i)]
\item $m$ divides an element in $\Omega$.
\item If $\phi(m)=\mbox{rank}(T(X))$, then $m\in\Omega$.
\item Conversely, for any $m\in\Omega$, there is a unique K3 surface as above with $\phi(m)=\mbox{rank}(T(X))$.
\end{enumerate}
\end{Theorem}

The non-unimodular case is less uniform. 
In order to formulate analogues of (ii) and (iii), we distinguish the following two sets
\[
\Omega_1=\{3,9,27,5,25,7,11,13,17,19\},\;\;\;\; \Omega_2=\{2,4,8,16\}.
\]
Here the uniqueness part of (iii) is due to Machida - Oguiso \cite{MO} for $m=25$ and Oguiso - Zhang \cite{OZ} for all other cases.

\begin{Theorem}[non-unimodular case]\label{Thm:non}
Let $X$ be an algebraic K3 surface with an automorphism $\eta$ of order $m>1$. Assume that $\NS(X)$ is non-unimodular and $\eta$ acts trivially on $\NS(X)$. Then
\begin{enumerate}[(i)]
\item $m\in\Omega_1\cup\Omega_2$.
\item If $\phi(m)=\mbox{rank}(T(X))$, then $m\in\Omega_1$.
\item Conversely, for any $m\in\Omega_1$, there is a unique K3 surface as above with $\phi(m)=\mbox{rank}(T(X))$.
\end{enumerate}
\end{Theorem}

The elements of $\Omega_2$ are missing in (ii) and (iii).
For these $2$-powers, the next rank of $T(X)$ compatible with (\ref{eq:phi}) is $\mbox{rank}(T(X))=m$.
Theorem~\ref{thm} gives a complete classification of this case.
The proof of Theorem~\ref{thm} is given in section \ref{s:proof}. 
First we recall some lattice theory.

\subsection{Discriminant group and $p$-elementary lattices}
\label{s:p}

Any integral lattice $L$ has a canonical embedding into its dual lattice $L^\vee$. We define the discriminant group $A_L$ of $L$ as the quotient
\[
A_L = L^\vee / L.
\]
If $L$ is non-degenerate, then $A_L$ is a finite abelian group. In the present situation, we consider a K3 surface $X$ with perpendicular lattices $\NS(X), T(X)$. Nikulin \cite{N} proved that
\begin{eqnarray}\label{eq:A}
A_{\NS(X)} \cong A_{T(X)}.
\end{eqnarray}
We say that a lattice $L$ is $p$-elementary (with $p$ prime) if $A_L$ is a $p$-elementary abelian group. 
The main step towards establishing Theorem~\ref{Thm:non} is the following result due to Vorontsov \cite{Vo}:

\begin{Theorem}\label{Thm:p-elem}
Let $X$ be an algebraic K3 surface with an automorphism $\eta$ of order $m>1$. Assume that $\NS(X)$ is non-unimodular and $\eta$ acts trivially on $\NS(X)$. Then
\begin{enumerate}[(i)]
\item
$m=p^k$ for some prime $p$.
\item
$\NS(X)$ is a $p$-elementary lattice.
\end{enumerate}
\end{Theorem}

This result readily puts us in a position to prove Theorem~\ref{thm} for $m=2$: 
Here $T(X)$ is positive-definite of rank two. Since it is $2$-elementary by Theorem~\ref{Thm:p-elem}, we obtain $T(X)=\langle 2\rangle^2$. By the Torelli theorem, this determines a unique complex K3 surface up to isomorphism. 

On the other hand, the given elliptic K3 surface has exactly four singular fibres, two each of types $I_2$ and $II^*$. 
Hence $\rho(X)=20$ and $\NS(X)$ is as claimed (cf.~Sect.~\ref{s:elliptic}). This implies that $T(X)=\langle 2\rangle^2$. 
Consider the elliptic involution $\eta: y\mapsto -y$, which acts trivially on the singular fibres. Since $\NS(X)$ is generated by fibre components and the zero section, $\eta$ operates trivially on $\NS(X)$. This completes the proof of the case $m=2$ of Theorem~\ref{thm}.

\subsection{$2$-elementary lattices}
\label{s:2}

By Theorem~\ref{Thm:p-elem}, we have to deal with $2$-elementary lattices to prove Theorem~\ref{thm}. These have been studied in great detail by Nikulin \cite{N-2}. To recall his classification result, we introduce the following notation.

For a non-degenerate integral lattice $L$, let $\ell(L)$ denote the minimal number of generators of the discriminant group $A_L$. 
Consider the induced quadratic form $\langle \cdot,\cdot\rangle$ on $L^\vee$.
If $L$ is $2$-elementary, define
\[
\delta(L) = 
\begin{cases}
0, & \text{if $\langle x,x\rangle\in\Z$ for all $x\in L^\vee$},\\
1, & \text{otherwise}.
\end{cases}
\]

\begin{Example}[Dynkin diagrams]
The $2$-elementary Dynkin diagrams, with their values of $\ell$ and $\delta$, are as follows:
$$
\begin{array}{c||c|c|c|c|c}
\hline
\mbox{type} & A_1 & E_7 & E_8 & D_{4n}\, (n>0) & D_{4n+2} \, (n>0)\\
\hline \hline
\ell & 1 & 1 & 0 & 2 & 2\\
\hline
\delta & 1 & 1 & 0 & 0 & 1\\
\hline
\end{array}
$$
\end{Example}

\begin{Theorem}[Nikulin {\cite[Theorem~4.3.2]{N-2}}]\label{Thm:2-elem}
Let $L$ be an even $2$-elementary lattice of rank $r$ and signature $(1, r-1)$. Then the isomorphism class of $L$ is determined by the triple $(r, \ell(L), \delta(L))$.
\end{Theorem}

In the same paragraph \cite[\S 4.3]{N-2}, Nikulin gives precise conditions for the existence of an even $2$-elementary lattice  $L$ with prescribed $(r, \ell(L), \delta(L))$. In our situation, we furthermore have to take into account that
\[
\ell(\NS(X))=\ell(T(X))
\]
by (\ref{eq:A}). In particular we obtain the trivial bound 
\begin{eqnarray}\label{eq:min}
\ell(\NS(X))\leq\min(\mbox{rank}(\NS(X)), \mbox{rank}(T(X))).
\end{eqnarray}
With this bound and Theorem~\ref{Thm:2-elem}, we can easily list all $2$-elementary lattices which could possibly be associated to the non-unimodular K3 surfaces in Theorem~\ref{thm}. In Table \ref{tab}, we only give the hypothetical N\'eron-Severi lattices. All other triples $(r, \ell, \delta)$ are ruled out by Nikulin's  statement in \cite[\S 4.3]{N-2} and (\ref{eq:min}).

\begin{table}
$$
\begin{array}{c|c|c}
\hline
m & (r, \ell, \delta) & L\\
\hline
\hline
2 & (20, 2, 1) & U+A_1^2+E_8^2\\
\hline
4 & (18, 2, 0) & U+D_8+E_8\\
& (18,2,1) & U+A_1+E_7+E_8\\
& (18,4,0) & U + D_8^2\\
& (18,4,1) & U+A_1^2+E_7^2\\
\hline
8 & (14,2,0) & U+D_4+E_8\\
& (14,4,0) & U+D_4+D_8\\
& (14,4,1) & U+A_1^4+E_8\\
& (14,6,0) & U+D_4^3\\
& (14,6,1) & U+A_1^4+D_8\\
& (14,8,0) & U(2) + D_4^3\\
& (14,8,1) & U+A_1^4+D_4^2\\
\hline
16 & (6,2,0) & U+D_4\\
& (6,4,0) & U(2)+D_4\\
& (6,4,1) & U+A_1^4\\
& (6,6,1) & U(2)+A_1^4\\
\hline
\end{array}
$$
\caption{The $2$-elementary lattices $L$ possibly equalling $\NS(X)$ for non-unimodular K3 surfaces $X$ with $m=\mbox{rank}(T(X))=2^k$.}
\label{tab}
\end{table}

The corresponding transcendental lattices are easily computed by comparing the discriminant forms, following the theory developed by Nikulin \cite{N}.
In all present cases, there is only one class of lattices per genus, so that the discriminant form determines the lattice up to isometry.
Thus we verify the transcendental lattices given in Theorem~\ref{thm}.
For the non-symplectic K3 surfaces from Theorem \ref{Thm:uni} and \ref{Thm:non}, analogous arguments are presented in detail in \cite{LSY}.

We chose to write the N\'eron-Severi lattices in a very particular way, always involving $U$ or $U(2)$. The reason for this will become clear in the next section when we turn to elliptic fibrations.

Let $X$ be a K3 surface with $\NS(X)$ $2$-elementary. Nikulin \cite{N-2} showed that $X$ admits an involution $\iota$ such that
\[
\iota^*|_{\NS(X)} = 1,\;\;\; \iota^*|_{T(X)} = -1.
\]
On the K3 surfaces from Theorem~\ref{thm}, we will consider $\iota=\eta^{m/2}$. Then we will study the fixed curve
\[
\Theta=\text{Fix}(\iota).
\]
By the  Torelli theorem, $\iota$ is unique. Hence $\text{Aut}(X)$ is the centraliser of $\iota$. In particular, $\text{Aut}(X)$ maps $\Theta$ onto itself, so the curve $\Theta$ will be fixed by $\eta$.

\begin{Theorem}[Nikulin {\cite[Theorem~4.2.2]{N-2}}]\label{Thm:C}
$\Theta$ is a nonsingular curve. It decomposes into disjoint components depending on the triple $(\mbox{rank}(\NS(X)), \ell(\NS(X)), \delta(\NS(X)))=(r, \ell, \delta)$:
\[
\Theta= 
\begin{cases}
\emptyset & \text{ if } (r, \ell, \delta) = (10,10,0)\\
C_1 + C_2 & \text{ if } (r, \ell, \delta) = (10,8,0)\\
C + \sum_{i=1}^n B_i & \text{ otherwise}.
\end{cases}
\]
Here $C_1, C_2$ are smooth curves of genus one.  $C$ denotes a smooth curve of genus $g=(22-r-\ell)/2$. The $B_i$ are smooth rational curves, $n=(r-\ell)/2$.
\end{Theorem}

\subsection{Elliptic fibrations}
\label{s:elliptic}

K3 surfaces can admit several elliptic fibrations onto $\PP^1$. Here we further have to distinguish whether a given fibration has a section. If so, we denote it by $O$. Then the general fibre $F$ is an elliptic curve with the intersection point $F\cap O$ as origin of the group law. On an elliptic K3 surface, $O^2=-2$. Hence $F$ and  $O$ generate the hyperbolic lattice $U$.

We want to formulate a converse statement so that from the N\'eron-Severi lattices in Table \ref{tab} we can deduce the existence of an elliptic fibration. For this we identify the reducible singular fibres with Dynkin diagrams. If there is a section $O$, the identification is achieved by omitting the fibre components that meet $O$. In general, one omits a simple component on each singular fibre (unless there are multiple fibres). Then one just draws the intersection graph. The following table pairs the type of the singular fibre in Kodaira's notation with the corresponding Dynkin diagram:
$$
\begin{array}{c||c|c|c|c|c|c|c}
\hline
\text{fibre type} & I_2, III & I_3, IV & I_n (n>3) & I_n^* (n\geq 0) & IV^* & III^* & II^*\\
\hline \hline
\text{Dynkin diagram} & A_1 & A_2 & A_{n-1} & D_{n+4} & E_6 & E_7 & E_8\\
\hline
\end{array}
$$

\begin{Lemma}\label{Lem:U}
Let $X$ be a K3 surface. Assume that $\NS(X)=U + \Gamma_1+\hdots+\Gamma_n$ where each $\Gamma_i$ denotes a Dynkin diagram. Then $X$ admits an elliptic fibration with section and singular fibres corresponding to the $\Gamma_i$.
\end{Lemma}

A proof of this lemma can be found in \cite[Lemma~2.1, 2.2]{Ko}. The lemma applies to most lattices in Table \ref{tab}. Kond\=o also gave a generalisation for the remaining lattices which include a summand of $U(2)$. Here we need the extra information that the lattice is $2$-elementary.

\begin{Lemma}\label{Lem:U(2)}
Let $X$ be a K3 surface. Assume that $\NS(X)=U(2) + \Gamma_1+\hdots+\Gamma_n$ where $\Gamma_i=A_1, E_7, E_8, D_{4n} (n\geq 1)$. Then $X$ admits an elliptic fibration with singular fibres corresponding to the $\Gamma_i$.
\end{Lemma}

\begin{figure}[ht!]
\setlength{\unitlength}{.45in}
\begin{picture}(10,2.6)(-0.1,-0.3)

%
%
%

\thinlines

\put(5.5,1.7){\line(1,0){1}}
\put(5.5,1.4){\line(1,0){1}}
\put(5.5,0.3){\line(1,0){1}}
\put(5.5,0.6){\line(1,0){1}}

\thicklines
\put(5.7,2){\line(0,-1){2}}
\put(6.8,1.8){\makebox(0,0)[l]{$I_0^*$}}

%
%

\thinlines

\put(0.4,1.6){\line(1,0){0.6}}
\put(0.4,0.4){\line(1,0){0.6}}

\put(0,1.8){\line(1,0){0.8}}

\put(0,1){\line(1,0){0.7}}
\put(0,0.2){\line(1,0){0.8}}

\thicklines
\put(0.2,2){\line(0,-1){2}}


\put(0.6,2){\line(0,-1){0.6}}
\put(0.6,0){\line(0,1){0.6}}

\put(1.3,1.8){\makebox(0,0)[l]{$III^*$}}

\thinlines
\put(3.75,0.4){\line(1,0){0.45}}

\put(3,1.8){\line(1,0){0.8}}

\put(3,1.1){\line(1,0){0.7}}
\put(3,0.2){\line(1,0){0.75}}

\put(3.4,0.6){\line(1,0){0.7}}

\thicklines
\put(3.2,2){\line(0,-1){2}}


\put(3.6,2){\line(0,-1){0.6}}

\put(3.6,0){\line(0,1){0.8}}

\put(4.15,0.65){\makebox(0,0)[l]{{\tiny 3}}}

\put(3.9,0.7){\line(0,-1){0.45}}

\put(4.1,1.8){\makebox(0,0)[l]{$II^*$}}

%
%

%
%
%

\thinlines
\put(8,1.8){\line(1,0){1}}
\put(8,1.6){\line(1,0){1}}
\put(9,0.2){\line(1,0){1}}
\put(9,0.4){\line(1,0){1}}

\put(8,1.2){\line(1,0){0.8}}

\thicklines
\put(8.2,2){\line(0,-1){1}}

\put(8.6,1.4){\line(0,-1){0.8}}
\put(9.8,0){\line(0,1){1}}

\put(9.3,1.8){\makebox(0,0)[l]{$I_{4n-4}^*\;(n>1)$}}
\put(9,0.7){\makebox(0,0)[l]{$\hdots$}}

\end{picture}
\caption{Action of $\iota$ on non-reduced singular fibres}
\label{Fig:action}
\end{figure}
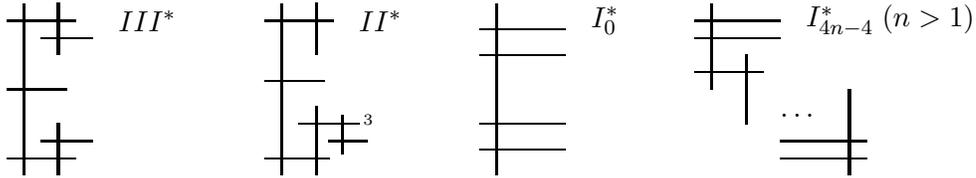

By the previous two lemmas, it suffices for our classification to consider \emph{elliptic} K3 surfaces with $\NS(X)$ $2$-elementary as in Table \ref{tab}. By Nikulin \cite[\S 4.2]{N-2}, any such K3 surface is equipped with an involution $\iota$ which acts trivially on $\NS(X)$. 
In particular, $\iota$ preserves the elliptic fibration and maps each section, if there is any, to itself. Kond\=o \cite[Lemma~2.3]{Ko} describes the operation on the singular fibres:

\begin{Lemma}\label{Lem:sing}
Let $X$ be an elliptic K3 surface with $\NS(X)$ $2$-elementary as in Lemma~\ref{Lem:U(2)}. 
\begin{enumerate}[(i)]
\item
The involution $\iota$ acts on the simple components of the singular fibres as an automorphism of order two. 
\item
On the multiple components, $\iota$ acts either as identity or as involution of order two.
The precise pattern is as follows:
$\iota$ acts as identity on the multiple components meeting simple components;
from there on, its action alternates between the two possibilities as depicted in Figure \ref{Fig:action}.
\end{enumerate}
\end{Lemma}

For the relevant non-reduced fibre types corresponding to $E_7, E_8, D_{4n} (n\geq 1)$, 
we sketch the action of $\iota$ on the fibre components in the following figure. 
Multiple components will be printed thick and vertically if $\iota$ acts as identity; 
all other components, in particular the simple ones, will appear horizontally in thin print.

By this lemma, we can identify many components of the fixed curve $\Theta$ of $\iota$ in Theorem~\ref{Thm:C} as fibre components (or as the section $O$ if $\iota|_{\PP^1}=1$). We will then argue using the remaining components. 

We conclude this section by recalling some basic facts from the theory of elliptic surfaces that will play a central role in our analysis:
\begin{itemize}
 \item 
By the Shioda-Tate formula, the N\'eron-Severi group $\NS(S)$ of an elliptic surface $S$ with section is generated by horizontal and vertical divisors, i.e.~sections and fibre components.
Hence permutations of reducible fibres induce a non-trivial action on $\NS(S)$.
\item
Over $\C$, there are exactly two elliptic curves that admit automorphisms of order $>2$ (respecting the group structure).
They have $j$-invariants $j=0, 1728$ and CM by the full integer rings in $\Q(\sqrt{-3})$ resp.~$\Q(\sqrt{-1})$.
\item
If the generic fibre of an elliptic surface with section has CM, then the fibration is isotrivial. 
This restricts the possible singular fibres.
For instance, if the generic fibre admits an automorphism of order four, i.e.~if $j=1728$, then the only possible singular fibre types are $III, I_0^*, III^*$.
\end{itemize}

\section{The classification result}
\label{s:proof}

This section gives the proof of Theorem \ref{thm}, divided into subsections corresponding to the four families of K3 surfaces occuring.
The single surface in case $m=2$ has already been treated in \ref{s:p}.

\subsection{Proof in the unimodular case}
\label{s:uni}

We are looking for all K3 surfaces $X$ with a non-symplectic automorphism $\eta$ under the following assumptions: $\eta$ acts trivially on $\NS(X)$, the order $m$ of $\eta$ is a $2$-power and $T_X$ has rank $m$ and is unimodular.

By Theorem~\ref{Thm:uni}, the only possibility is $m=\mbox{rank}(T(X))=4$. Hence by the classification of even unimodular lattices of given signature,
\[
\NS(X) = U + E_8^2,\;\;\; T(X) = U^2.
\]
By Lemma~\ref{Lem:U}, $X$ admits an elliptic fibration with section and two singular fibres of type $II^*$. 
Such K3 surfaces have been studied in great detail by Shioda in \cite{Sandwich}.
They are given by the Weierstrass equation
\begin{eqnarray}\label{eq:4-uni}
X:\;\;\; y^2 = x^3 - 3\,\lambda\,t^4\,x+ t^7 + \mu\,t^6+ t^5.
\end{eqnarray}
with $II^*$ fibres at $t=0$ and $t=\infty$. 
In particular, any such K3 surface $X$ admits a Shioda-Inose structure:
the quadratic base change $t\mapsto t^2$ results in another elliptic K3 surface. By \cite{Sandwich}, this is the Kummer surface of the product of two elliptic curves $E, E'$. 
Based on an argument by Inose \cite{Inose}, the elliptic curves are determined by their $j$-invariants through the parameters $\lambda, \mu$ (cf.~(\ref{eq:OO}) for  the special case $\mu=0$).

Because of the singular fibres of type $II^*$, the general fibre of $X$ does not admit an automorphism of order four. 
Hence $\eta$ has to operate non-trivially on the base curve $\PP^1$. Since $\eta$ preserves the elliptic fibration, we deduce $\mu=0$. This reduces (\ref{eq:4-uni}) to the equation in Theorem~\ref{thm}. Here $\eta$ is given by
\[
\eta:\;\; t\mapsto -t,\;\; x\mapsto -x,\;\; y\mapsto \sqrt{-1}\, y.
\]
For the corresponding elliptic curves, this implies that one of them, say $E$,  has $j(E)=1728$. 
Thus $E$ admits an automorphism of order four (which induces $\eta$). 

Finally we have to make sure that $\eta$ acts trivially on $\NS(X)$. This certainly holds true for $O, F$ and the two singular fibres of type $II^*$. However, we could have $\rho(X)>18$ (so that by (\ref{eq:phi}) already $\rho(X)=20$). In all such cases, one can see that there are some additional cycles that are not $\eta$-invariant. 

For instance, if $\lambda^3=1$, there are further reducible singular fibres of type $I_2$ at $t=\pm 1$. Hence $\eta$ interchanges them. In fact, the resulting surface is isomorphic to the K3 surface for $m=2$ in Theorem~\ref{thm}. The non-symplectic automorphism $\eta^2$ acts trivially on $\NS(X)$. In terms of the Shioda-Inose structure, this is exactly the case $E\cong E'$.

On the other hand, there could be additional sections. By the Shioda-Inose structure, this happens if and only if the two elliptic curves are isogenous, but not isomorphic. Since $E$ has CM and $j(E)=1728$, this is equivalent to $E'$ having complex multiplication by some order in $\Q(\sqrt{-1})$ ($\neq \Z[\sqrt{-1}]$). In terms of the parameter $\lambda$, both degenerate cases together can be expressed as follows:
\begin{eqnarray}\label{eq:OO}
\rho(X)=20 \Longleftrightarrow 1728\,\lambda^3 =j(\OO)\;\; \text{ for some order }\; \OO\subseteq\Z[\sqrt{-1}].
\end{eqnarray}

\subsection{Non-unimodular case $m=4$}
\label{s:non-4}

We first rule out all $2$-elementary lattices $L$ from Table~\ref{tab} except for one. 
Then we derive the family of K3 surfaces given in Theorem~\ref{thm} for the remaining lattice. 
In each hypothetical case, the assumption $\NS(X)=L$ guarantees an elliptic fibration on $X$ with certain singular fibres by Lemma~\ref{Lem:U}. We will always work with this fibration.

If $\NS(X)=U+A_1+E_7+E_8$ or $U+A_1^2+E_7^2$, then there are more than two reducible singular fibres. 
As explained, interchanging them induces a non-trivial action on $\NS(X)$.
Hence $\eta$ has at least three fixed points on the base curve $\PP^1$, so it operates trivially. 
Thus $\eta$ also fixes $O$.
Hence the general fibre is an elliptic curve with an automorphism of order four. I.e.~the fibration is isotrivial with $j=1728$. 

With $\NS(X)$ of the given shape, isotriviality is only possible in the second case with singular fibres $III, III^*$ twice each. By a M\"obius transformation, we move the singular fibres to $0,1,\gamma,\infty$. Then it follows from Tate's algorithm \cite{Tate} that the elliptic fibration is given up to isomorphism as
\begin{eqnarray}\label{eq:iso}
X:\;\;\; y^2 = x^3 + t^3\,(t-1)\,(t-\gamma)\,x.
\end{eqnarray}
The automorphism of order four is indeed operating trivially on the singular fibres. However, there is a two-torsion section $(0,0)$. 
Hence the lattice $U+A_1^2+E_7^2$ has index two in $\NS(X)$.  
At the end of this section, we will verify that $\NS(X)=U+D_8+E_8$. 


The final case to be ruled out is $\NS(X)=U+D_8^2$. Here we could again argue with an explicit Weierstrass equation. However, we decided to give a geometric proof that no such elliptic surface admits an automorphism of order four with trivial action on $\NS(X)$. The proof follows the lines of Kond\=o's arguments in \cite{Ko}.

We will use that $\eta$ acts non-trivially on $\PP^1$. Otherwise, the general fibre would have CM again which is not possible with singular fibres of type $I_4^*$. We let $\iota=\eta^2$. By Theorem~\ref{Thm:C}, there are smooth rational curves $B_i (i=1,\hdots,8)$ such that
\[
\mbox{Fix}(\iota) = \sum_{i=1}^8 B_i
\]
By Lemma~\ref{Lem:sing}, we may assume that $B_1=O$, and $B_2,\hdots,B_7$ are disjoint double components of the singular fibres. Moreover, there are eight isolated fixed points of $\iota$, one on each simple component of the singular fibres. Exactly two of these points lie on $O$. Hence the remaining six lie on $B=B_8$. 
Since $\eta$ operates trivially on $\NS(X)$ by assumption, each of these fixed points of $\iota=\eta^2$ is already a fixed point of $\eta$. 

We deduce that $B$ intersects the general fibre in three points. Since $\eta|_{\PP^1}\neq 1$, this implies $\eta|_{B}\neq 1$. Hence we can apply the Hurwitz formula to $B$ and $\eta$. With $d=\mbox{ord}(\eta|_B)$ it reads
\[
-2 = 2\, (g(B)-2) = d\, (2\,g(B/\eta) -2) + 6\,(d-1)\geq 4\,d-6.
\]
Since $d>1$, this gives the required contradiction. 

For the remaining lattice $\NS(X)=U+D_8+E_8$, we shall now derive the family of elliptic surfaces given in Theorem~\ref{thm}. Then we will check the compatibility with the isotrivial fibration (\ref{eq:iso}).

With singular fibres of type $I_4^*, II^*$ at $0, \infty$, Tate's algorithm predicts the following Weierstrass equation:
\begin{eqnarray}\label{eq:non-4}
X:\;\;\; y^2 = x^3 + (\mu\,t+\nu)\,t\,x^2 + \lambda\,t^4\,x+\gamma\,t^7,\;\;\;\nu\gamma\neq 0.
\end{eqnarray}
After scaling, we can assume $\nu=\gamma=1$.
Here $\eta|_{\PP^1}\neq 1$ for the same reason as before. Since $\eta$ preserves the elliptic fibration, we deduce $\mu=0$. Hence (\ref{eq:non-4}) reduces to the equation in Theorem~\ref{thm}. Then $\eta$ can be given as
\[
\eta:\;\; t\mapsto -t,\; x\mapsto -x,\; y\mapsto \sqrt{-1}\,y.
\]
This elliptic surface has discriminant
\[
\Delta =16\, t^{10}\,(27\,t^4-2\,\lambda\,(2\,\lambda^2+9)\,t^2-\lambda^2+4),
\]
so in general there are four singular fibres of type $I_1$. They degenerate exactly in the following cases: If $\lambda^2=3$, then there are two fibres of type $II$ instead. If $\lambda=\pm 2$, then two $I_1$'s merge with the fibre $I_4^*$ at $t=0$ to form $I_6^*$. Hence $\rho=20$, and we obtain the surface from case $m=2$. Note that $\eta$ does not act trivially on the fibre of type $I_6^*$ any more, but, of course, $\eta^2$ does.

We still have to show $\eta|_{\NS(X)}=1$ for a general choice of $\lambda$. For this, it suffices to verify that $\rho(X)=18$ (so that in particular $\NS(X)=U+D_8+E_8$, since a $II^*$ fibre does not admit any torsion sections). We show this using the smooth specialisation $X_0$ at $\lambda=0$. Clearly $18\leq\rho(X)\leq\rho(X_0)$. On $X_0$, we can take the square root of $\eta$: fixing a primitive eighth root of unity $\zeta$, we have
\[
\Aut(X)\ni\sqrt\eta:\;\; t\mapsto\zeta^2\,t,\; x\mapsto\zeta^2\,x,\; y\mapsto\zeta^3\,y.
\]
By (\ref{eq:phi}), $T(X_0)$ has rank at least four. I.e.~$\rho(X_0)\leq 18$, which implies the equality $\rho(X)=\rho(X_0)=18$. In particular $\eta$ operates trivially on $\NS(X)$ for general $\lambda$. This completes the proof of Theorem~\ref{thm} in case $m=4$. 

We conclude this section by checking the compatibility of the two elliptic fibrations (\ref{eq:iso}) and (\ref{eq:non-4}). We exhibit an alternative elliptic fibration on the K3 surfaces (\ref{eq:non-4}). For this we consider the affine chart $x=t^3\,u, y=t^3\,v$ of the triple blow-up at $(0,0,0)$:
\[
X:\;\;\; v^2 = t^3\,u^3 + t\,u^2 + \lambda\,t\,u+t.
\]
We now choose $u$ as a section. A simple variable change produces the Weierstrass equation
\[
X:\;\;\; v^2 = t^3 + u^3\,(u^2+\lambda\,u+1)\,t.
\]
This reveals the relation to the family of isotrivial elliptic fibrations (\ref{eq:iso}):
\begin{eqnarray}\label{eq:gamma}
\lambda = - \dfrac{1+\gamma}{\sqrt{\gamma}}.
\end{eqnarray}
In section \ref{s:arith}, we will use this relation to determine the K3 surfaces in the family with $\rho=20$ (cf.~Corollary~\ref{cor}). Those surfaces are excluded in Theorem~\ref{thm}.

\subsection{Proof of case $m=8$}
\label{s:8}

By the same methods as before, we can rule out the four cases $\NS(X)=L$ where $L$ can be written as sum of $U$ and at least three Dynkin diagrams: Here $\eta|_{\PP^1}=1$ and $\eta$ fixes $O$. Hence the general fibre is an elliptic curve with an automorphism of order 8, contradiction.

If $\NS(X)=U(2)+D_4^3$, we still have $\eta|_{\PP^1}=1$, but no section. Instead we work with the fixed curve of $\iota=\eta^4$. By Theorem~\ref{Thm:C}, there are disjoint smooth rational curves $B_1,\hdots,B_4$ such that
\[
\mbox{Fix}(\iota) = B_1 + \hdots + B_4.
\]
By Lemma~\ref{Lem:sing}, three of these curves are the multiple components of the $I_0^*$ fibres. Denote the remaining rational curve by $B$. Since $\eta$ operates trivially on $\NS(X)$, it fixes each single $B_i$. Moreover, $\eta$ has 12 fixed points where $B$ intersects the simple components of the $I_0^*$ fibres again by Lemma~\ref{Lem:sing}. We distinguish two cases depending on $d=\mbox{ord}(\eta|_B)$. 

If $d=1$, then the intersection with $B$ equips each fibre $F$ with four rational points fixed by $\eta$. Hence  $F$ is an elliptic curve  with an automorphism of order 8. As above this gives a contradiction. 

If $d>1$, then we establish a contradiction with the Hurwitz formula applied to $B$ with $\eta$ and the 12 fixed points.

To rule out $\NS(X)=U+D_4+D_8$, we apply a similar argument. Here $\eta|_{\PP^1}$ has order at least four, since otherwise the general fibre would have an automorphism $\eta^2$ of order four. By Theorem~\ref{Thm:C} there is a smooth curve $C$ of genus two such that with disjoint multiple components $B_1,\hdots,B_4$ of the singular fibres
\[
\mbox{Fix}(\eta^4) = C + O + B_1 +\hdots + B_4.
\]
Again $C$ is fixed by $\eta$. Since $\eta$ acts trivially on $\NS(X)$, it has six fixed points where $C$ intersects the simple components of the singular fibres. In particular $C.F>0$, so that $d=\mbox{ord}(\eta|_C)\geq\mbox{ord}(\eta|_{\PP^1})\geq 4$. Now we apply the Hurwitz formula
\[
2=2\,g(C)-2 = d\,(2\,g(C/\eta)-2) + 6\,(d-1)\geq 4\,d-6
\]
to establish the contradiction $d\leq 2$.

Finally we derive the Weierstrass form for the family of elliptic surfaces with the remaining lattice $\NS(X)=U+D_4+E_8$. We locate the singular fibres at $0, \infty$.  Then the Tate algorithm predicts the Weierstrass equation
\[
X:\;\;\; y^2 = x^3 + A(t)\,t\,x^2 + B(t)\,t^2\,x+C(t)\,t^3
\]
with $\deg(A(t))\leq 1, \deg(B(t))\leq 2, \deg(C(t))=4$. 
After translating $x$, we can assume $C(0)=0$. Then the fibre at $t=0$ has type $I_0^*$ if and only if $B(0)\neq 0$.

Now we use that by the same arguments as before $\mbox{ord}(\eta|_{\PP^1})\geq 4$. 
Since $\eta$ preserves the elliptic fibration, it acts as multiplication by some scalar on the polynomials $A, B, C$.
From the low degrees (and $C(0)=0$), it follows that $A, B, C$ are all monomials. 
By the above conditions, we may assume that, after scaling,
\[
C(t)=t^4,\;\; B(t)=1.
\]
Here $\eta$ has to operate as $t\mapsto \sqrt{-1}\,t, x\mapsto\sqrt{-1}\,x$. Hence $A(t)=\lambda$, giving the equation from Theorem~\ref{thm}. The discriminant is
\[
\Delta =16\, t^{6}\,(27\,t^8-2\,\lambda\,(2\,\lambda^2+9)\,t^4-\lambda^2+4),
\]
so in general there are 8 fibres of type $I_1$. Degenerations occur exactly at $\lambda=\pm 2$ where four of them collapse with $I_0^*$ to form $I_4^*$, and at $\lambda^2=3$ with four $II$'s instead. 
Hence for $\lambda$ general, $\eta$ acts trivially on the reducible singular fibres. For the remaining claims about the general member $X$, the same argument with the smooth specialisation $X_0$ at $\lambda=0$ applies as in section \ref{s:non-4}.

\subsection{Proof of case $m=16$}
\label{s:16}

If $\NS(X)=U+A_1^4$ or $U(2)+A_1^4$, we again have $\eta|_{\PP^1}=1$. In the first case, there is a section (fixed by $\eta$). Hence the general fibre is an elliptic curve with an automorphism of order 16, contradiction. 

In the second case, $C=\mbox{Fix}(\eta^8)$ is a smooth curve of genus five by Theorem~\ref{Thm:C}. By Lemma~\ref{Lem:sing}, $C$ meets each component of the reducible singular fibres in two points. Hence
$C$ intersects the general fibre $F$ in four points.  In particular, these are fixed by $\eta^4$. This means that $F$ is an elliptic curve with an automorphism $\eta^4|_F$ of order four which fixes four points. This is impossible. 

We now consider the lattice $\NS(X)=U(2)+D_4$. By Theorem~\ref{Thm:C}, 
\[
\mbox{Fix}(\eta^8) = B + C
\]
where $B$ is the multiple component of the special fibre and $C$ is a smooth curve of genus 6. By Lemma~\ref{Lem:sing}, $C$ meets each simple component of the special fibre in a point which is actually fixed by $\eta$. Hence $C.F=4$ and $\#\mbox{Fix}(\eta)\geq 4$. 

This implies that the order of the $\eta$-action on $\PP^1$ is at most four. Otherwise $d=\mbox{ord}(\eta|_C)\geq\mbox{ord}(\eta|_{\PP^1})\geq 8$, since $C.F>0$. Then the Hurwitz formula would give
\[
10 = 2\,g(C)-2 = d\,(2\,g(C/\eta)-2) + 4 (d-1) \geq 4\,d-6,
\]
so $d\leq 4$, contradiction. 

Define $\xi=\eta^4$ with trivial action on $\PP^1$. If $\xi|_C=1$, then the general fibre is an elliptic curve with an automorphism of order four, fixing four points. As above, this gives a contradiction. 

If $\mbox{ord}(\xi|_C)=d>1$, then we apply the Hurwitz formula to $C$ and $\xi$. Here we use that $\xi$ has at least 22 fixed points on $C$: 
the nodes of the singular fibres of type $I_1$ or the cusps and one further point on the $II$ fibres plus the four intersection points with the $I_0^*$ fibre. Hence the Hurwitz formula
\[
10 = 2\,g(C)-2 = d\,(2\,(g(C/\xi)-2) + 22 (d-1) \geq 20\,d-22
\]
gives the contradiction $d<2$. This completes the non-existence proof.

It remains to derive the family of elliptic surfaces with $\NS(X)=U+D_4$ from Theorem~\ref{thm}. We work with an elliptic fibration where we locate the special fibre at $t=0$. By the same arguments as before, $\eta|_{\PP^1}$ has order at least 8. Hence there are 8 singular fibres of type $II$ or 16 $I_1$ which are interchanged by $\eta$. Since $e(X)=\sum_F e(F)=24$, there remains one singular fibre of type $II$ which is fixed by $\eta$. We locate it at $\infty$ with cusp at the origin.
Then Tate's algorithm gives
\begin{eqnarray}\label{eq:16}
X:\;\;\; y^2 = x^3 + A(t)\,t\,x^2 + B(t)\,t^2\,x+C(t)\,t^3
\end{eqnarray}
with $\deg(A(t))\leq 2,\, \deg(B(t))\leq 5,\, \deg(C(t))=8$. After translating $x$, we can assume $C(0)=0$. Then the fibre has type $I_0^*$ at $t=0$ if and only if $B(0)\neq 0$. 

As before, $\eta$ acts as some scalar multiplication on the polynomials $A, B, C$.
Since $\eta$ has order at least $8$ on $\PP^1$, we deduce that each polynomial is in fact a monomial due to its small degree. After normalising, we obtain
\[
C(t) = t^{8},\;\; B(t)=1,\;\; A(t)=\lambda.
\]
Thus (\ref{eq:16}) reduces to the claimed family of elliptic K3 surfaces. 
The discriminant is
\[
\Delta = 16\, t^{6}\,(27\,t^{16}-2\,\lambda\,(2\,\lambda^2+9)\,t^8-\lambda^2+4),
\]
so in general there are 16 fibres of type $I_1$. Degenerations occur exactly in the two usual cases: at $\lambda=\pm 2$ where eight $I_1$'s merge with $I_0^*$ to constitute $I_8^*$, and at $\lambda^2=3$ with eight $II$'s instead. 

The non-symplectic automorphism $\eta$ involves a primitive $16$-th root of unity $\zeta$:
\[
\eta:\;\;\;x\mapsto \zeta^2\,x, \; t\mapsto \zeta^2\,t,\; y\mapsto\zeta^3\,y.
\]
For $\lambda$ general, $\eta$ acts trivially on the reducible singular fibres. All other claims about the general member $X$ can be proved with the smooth specialisation $X_0$ at $\lambda=0$ as in \ref{s:non-4} and \ref{s:8}.

\section{Arithmetic aspects}
\label{s:arith}

In this section we will discuss arithmetic aspects of the K3 surfaces in Theorem~\ref{thm}. 
In particular, we will show that each family contains at least three members of CM type. 

First we note that the surface for $m=2$ in Theorem~\ref{thm} has $\rho=20$, hence is modular by \cite{L}. The associated Hecke eigenform has weight 3 and level 16 as given in \cite[Tab.~1]{S-CM}. 

In all other cases of Theorem~\ref{thm}, we are concerned with one-dimensional families of K3 surfaces. Hence any relation to automorphic forms (as predicted by the Langlands program) will be more complicated. The transcendental lattice gives rise to a compatible system of $m$-dimensional Galois representations $\varrho$ over $\Q$. 
However, we can still reduce to two-dimensional Galois representations over some extension of $\Q$.
For this we fix a primitive root of unity $\zeta_m$ of order $m$. 

\begin{Proposition}\label{Prop:cyclo}
Let $X$ be a K3 surface over a number field $K$ with a non-symplectic automorphism of order $m$. Then the Galois representation $\varrho$ associated to $T(X)$ splits into $m$ equidimensional Galois representations over $K(\zeta_m)$.
\end{Proposition}

The proposition relies on the fact that the non-symplectic automorphism endows $T(X)$ with the structure of a  $\Z[\zeta_m]$-module (leading to (\ref{eq:phi})). This property carries over to the Galois representations. 
Applied to the families from Theorem~\ref{thm}, Proposition~\ref{Prop:cyclo} produces two-dimensional Galois representations over $\Q(\zeta_m)$.

In the unimodular case of $m=4$, we can describe the two-dimensional Galois representation explicitly. From the Shioda-Inose structure with the elliptic curves $E, E'$, it follows that 
\begin{eqnarray}\label{eq:T-4}
T(X) = H^1(E) \otimes H^1(E')
\end{eqnarray}
if $\rho(X)=18$ (cf.~(\ref{eq:OO})). Over some extension, this relation translates into Galois representations. Since $E$ has CM, there is a Hecke character $\psi$ over $\Q(\sqrt{-1})$ associated. Then $H^1(E)=\mbox{Ind}_\Q^{\Q(\sqrt{-1})} \psi$. Hence $\varrho$ is induced by $\psi\otimes H^1(E')$.

\begin{Definition}
Let $X$ be a smooth projective surface over a number field $K$. We say that $X$ has \emph{CM type} if over some finite extension of $K$ the Galois representation $\varrho$ associated to $T(X)$ splits into one-dimensional Galois representations.
\end{Definition}

By (\ref{eq:OO}) and (\ref{eq:T-4}), a member of the unimodular family from Theorem~\ref{thm} is of CM type if and only if the elliptic curve $E'$ has CM as well. Here the Picard number jumps to 20 if and only if the CM field is $\Q(\sqrt{-1})$. The specialisation  $X_0$ at $\lambda=0$ of CM type has been studied in \cite{LSY}. Here $j(E')=0$ by (\ref{eq:OO}). Hence $E'$ admits an automorphism of order three. Together with $\eta$, this induces a non-symplectic automorphism of order $12$ on $X_0$.

Our next aim is to investigate CM type surfaces in the other families from Theorem~\ref{thm}. We start with the non-unimodular family for $m=4$. Thanks to the relation (\ref{eq:gamma}), we can work with the model $X_\gamma$ from (\ref{eq:iso}). We want to establish a structure similar to (\ref{eq:T-4}). Here we use that $X_\gamma$ is an isotrivial elliptic surface with smooth fibre $E$ of $j(E)=1728$.

For $\gamma\neq 0,1$, we apply the following base change to the elliptic surface $X_\gamma\to\PP^1$:
\begin{eqnarray*}
C_\gamma = \{v^4=u\,(u-1)\,(u-\gamma)\} & \to & \PP^1\\
(u,v) & \mapsto & u
\end{eqnarray*}
The base change results in the product $E\times C_\gamma$. This induces an embedding of $T(X_\gamma)$ into $H^1(E)\otimes H^1(C_\gamma)$. The involution $v\mapsto -v$ identifies a summand of $H^1(C_\gamma)$ coming from the elliptic curve
\[
E_\gamma: \;\;\; w^2 = u\,(u-1)\,(u-\gamma).
\]
Moreover there are three involutions that permute the points above $0, 1, \gamma, \infty$ pairwise (e.g.~$(u, w) \mapsto (\gamma/u, w/u)$).
For each involution, the quotient is an elliptic curve with an automorphism of order four, i.e.~it is isomorphic to $E$.
These quotients provide enough information to conclude that $\mbox{Jac}(C_\gamma)$ is isogenous to the product $E_\gamma\times E^2$ (an argument sketched to us once by R.~Kloosterman). 
Since $T(X_\gamma)$ has rank four  in general, but $T(E\times E)$ has rank two, we obtain the following structure:

\begin{Lemma}\label{lem}
The general surface $X_\gamma$ has \;\;$T(X_\gamma) = H^1(E) \otimes H^1(E_\gamma)$.
\end{Lemma}

In the above construction, we had to exclude $\gamma=1$ which corresponds to $\lambda=\pm 2$. As we know, that specialisation agrees with the surface for $m=2$ from Theorem~\ref{thm}. 

\begin{Corollary}\label{cor}
\begin{enumerate}[(i)]
\item 
The surface $X_\gamma$ has CM type if and only if $\gamma=1$ or $E_\gamma$ has CM. 
\item
$\rho(X_\gamma)=20$ if and only if  $\gamma=1$ or $E_\gamma$ has CM in $\Q(\sqrt{-1})$.
\end{enumerate}
\end{Corollary}

We shall now study the other non-unimodular families. Here we investigate the special members from the previous sections: $X_0$ at $\lambda=0$ and the degenerations at $\lambda^2=3,4$. 
Each surface can be shown to have CM type using Proposition~\ref{Prop:cyclo}. 
Below we will give an alternative proof by exhibiting a covering by a  Fermat surface.  This will also enable us to determine the zeta function.

%

\begin{Remark}
In each non-unimodular family, the surfaces at $\lambda$ and $-\lambda$ are isomorphic via $\sqrt\eta$. Hence it makes sense to refer to the specialisations $X_2$ and  $X_{\sqrt{3}}$ in the following. In fact, there are models of the families in terms of $\mu=\lambda^2$. For instance, one obtains for $m=4$
\[
X:\;\;\; y^2 = x^3 + \frac 1\mu \,t\,x^2 + t^4 x+t^7.
\]
In this model, the member at $\mu=0$ degenerates. Therefore we decided to use the given models with symmetry $\lambda\leftrightarrow-\lambda$.
\end{Remark}

\begin{Lemma}\label{Lem:CM}
In the non-unimodular families, $X_0, X_2$ and $X_{\sqrt{3}}$ have CM type.
\end{Lemma}

\emph{Proof:}
We have seen that $X_0$ admits a non-symplectic automorphism $\sqrt\eta$ of order $2m$. Hence the claim follows from Proposition~\ref{Prop:cyclo}. On $X_2$, the singular fibres degenerate in such a way that $T(X_4)$ has only rank $m/2$. Hence the same proposition applies.

For $m=4$, the surface $X_{\sqrt{3}}$ has CM type by Corollary~\ref{cor}. By (\ref{eq:gamma}), $\lambda=\sqrt{3}$ corresponds to $\gamma$ being a primitive sixth root of unity. Hence $E_\gamma$ has CM with $j=0$, since in general 
\[
j(E_\gamma)=2^8\,\dfrac{(\gamma^2-\gamma+1)^3}{\gamma^2\,(\gamma-1)^2}.
\]

On the other hand, the elliptic surfaces $X_{\sqrt{3}}$ for $m=8, 16$ are isotrivial with $j=0$. Hence the general fibre admits an automorphism $\omega$ of order three. As $\omega$ is non-symplectic, $\eta\omega$ has order $3m$. 
This implies the claim by Proposition~\ref{Prop:cyclo}.

Here isotriviality is a consequence of the number of fibres of types $II$ and  $II^*$. After completing the cube so that the coefficient of $x^2$ vanishes, the coefficient $B(t)$ of $x$ has total multiplicity 9 at the singular fibres. Since the $\deg(B(t))\leq 8$, $B\equiv 0$. Up to scaling, we obtain the Weierstrass equation
\begin{eqnarray}\label{eq:3}
X_{\sqrt{3}}:\;\;\; y^2 = x^3 +\sqrt\Delta \;\;\;\;\; (m=8, 16).
\end{eqnarray}

\subsection{Fermat surfaces}
\label{s:Fermat}

The prototype surfaces of CM type are Fermat surfaces. Here the action of roots of unity on coordinates provides a motivic decomposition of $H^2$ into one-dimensional eigenspaces. 
Following Weil \cite{Weil}, these eigenspaces correspond to Jacobi sums. 
Shioda \cite{Sh} showed that these properties carry over to Delsarte surfaces, i.e.~surfaces in $\PP^3$ defined by a polynomial with four terms. In the next section we will show that all surfaces in Lemma~\ref{Lem:CM} are Delsarte surfaces. Then we determine their zeta functions over finite fields. Since our arguments follow the same lines as \cite{LSY}, we will omit the details. 
Most of these ideas go back to N.~Katz, Ogus and Weil.


Let $S_n$ denote the complex Fermat surface of degree $n$:
\begin{eqnarray*}\label{eq:S_n}
S_n:\;\;\; \{x_0^n+x_1^n+x_2^n+x_3^n=0\}\subset\PP^3
\end{eqnarray*}
For $n> 4$, $S_n$ has general type while $S_4$ is a K3 surface with $\rho=20$. The $n$-th roots of unity act on coordinates as $\mu_n^3$. This induces a decomposition of $H^2(S_n)$ into one-dimensional eigenspaces $V(\alpha)$ with character. Here $\alpha$ runs through the character group
\[
\A_n:=
\{\alpha=(a_0,a_1,a_2,a_3)\in(\Z/n\Z)^4\,|\, a_i\not\equiv 0\pmod n,\,\sum_{i=0}^3 a_i\equiv 0\pmod n\,\}.
\]
Let $(\Z/n\Z)^*$ operate on $\A_n$ coordinatewise by multiplication. Let $\TT_n\subset\A_n$ consist of all those $\alpha\in\A_n$ such that the $(\Z/n\Z)^*$-orbit of $\alpha$ contains an element $(b_0,\hdots,b_3)$ with canonical representatives $0<b_i<n$ and
\[
\sum_{i=0}^3 b_i\neq 2\,n.
\]
Then the eigenspace $V(\alpha)$ is transcendental if and only if $\alpha\in\TT_n$. We obtain
\[
T(S_n) = \bigoplus_{\alpha\in\TT_n} V(\alpha).
\]
Weil \cite{Weil} showed that these eigenspaces correspond to Hecke characters over $\Q(\zeta_n)$. These can be expressed in terms of Jacobi sums. Given a  prime $p\nmid n$, choose $q=p^r\equiv 1\mod n$, so that there is a primitive character 
\[
\chi: \F_q^*\to \C^*
\]
of order $n$. For $\alpha\in\mathfrak{A}_m$, define the Jacobi sum
\[
\jj(\alpha) = \sum_\text{\small $\begin{matrix} v_1, v_2, v_3\in\F_q^*\\ v_1+v_2+v_3=-1\end{matrix}$} \chi(v_1)^{a_1} \chi(v_2)^{a_2}\chi(v_3)^{a_3}.
\]
\begin{Theorem}[Weil]
\label{Thm:Weil}
The Fermat surface $S_n$ over $\F_q$ has the following zeta function:
\[
\zeta(S_n/\F_q, T) = \dfrac 1{(1-T)\, P(T)\, (1-q^2\, T)}
\]
where
\[
P(T) = (1-q\,T)  \prod_{\alpha\in\mathfrak{A}_n} (1-\jj(\alpha)\,T).
\]
\end{Theorem}

\subsection{Zeta functions}
\label{s:zeta}

Shioda \cite{Sh} showed that the motivic decomposition of Fermat surfaces carries over to Delsarte surfaces, i.e.~surfaces in $\PP^3$ defined by a polynomial with four terms. Here we apply these ideas to the K3 surfaces in Lemma~\ref{Lem:CM} and  determine their zeta functions.

Here we will not consider $X_2$ or $X_{\sqrt{3}}$ for $m=4$. The former has $\rho=20$ and thus equals the surface for $m=2$. Hence the essential factor of the zeta function is given by the newform of weight 3 and level 16. The zeta function of $X_{\sqrt{3}}$ can be obtained from Lemma~\ref{lem} through $E_\gamma$ with $j(E_\gamma)=0$.

\begin{Lemma}\label{Lem:cover}
Consider the specialisations $X_0, X_2, X_{\sqrt{3}}$ in the non-unimodular families.
Except for $X_2, X_{\sqrt{3}}$ in case $m=4$, 
each surface  is covered by a Fermat surface.
\end{Lemma}

We first show that the surfaces are Delsarte surfaces. This implies the claim by \cite{Sh}, but we will also give the explicit covering maps.

The surfaces $X_0$ are visibly Delsarte surfaces. Now we let $m=8, 16$. On $X_{-2}$, the translation $x\mapsto x-t$ produces the representation as a Delsarte surface:
\begin{eqnarray}\label{eq:D-2}
X_{-2}:\;\;\; y^2 = x^3 + t\,x^2 + t^{3+m/2}.
\end{eqnarray}
Note that for $m=8$, this produces exactly $X_0$ from the $m=4$ case. 
On $X_{\sqrt{3}}$, the elliptic fibration (\ref{eq:3}) is a Delsarte model. After a variable change over $\Q(\zeta_m, 3^{1/m})$, the fibration becomes
\begin{eqnarray}\label{eq:D-3}
X_{\sqrt{3}}:\;\;\; y^2 = x^3 + t^3 + t^{3+m/2}.
\end{eqnarray}
For the covering maps, we will always work in the following affine chart of $S_n$:
\begin{eqnarray}\label{eq:aff}
S_n:\;\;\; u^n+v^n+w^n+1=0
\end{eqnarray}
For the Delsarte surfaces $X_2, X_{\sqrt{3}}$, we will employ the above affine models (\ref{eq:D-2}), (\ref{eq:D-3}). Then we write $y, -x, -t$ as functions of $u, v, w$. 
$$
\begin{array}{c|c||c||c|c|c}
\hline
m & \lambda & n & y & -x & -t\\
\hline \hline
4 & 0 & 8 & u^4\,v^{14}/w^{21} & v^{12}/w^{14} & v^4/w^6\\
\hline
8 & 0 & 16 & u^8\,v^7/w^{21} & v^{10}/w^{14} & v^2/w^6\\
 & \sqrt{3} & 24 & u^{12}\,w^9 & v^8\,w^6 & w^6\\
\hline
16 & 0 & 32 & u^{16}\,v^{11}/w^{33} & v^{18}/w^{22} & v^2/w^6\\
& -2 & 16 & u^8\,v^{22}/w^{33} & v^{20}/w^{22} & v^4/w^6\\
& \sqrt{3} & 48 & u^{24}\,w^9 & v^{16}\,w^6 & w^6\\
\hline
\end{array}
$$
In each of the above cases, let $G$ denote the subgroup of $\mu_n^3$ which leaves the coordinates $y, x, t$ invariant. 
It follows that the Delsarte surface $X$ is birationally given as the quotient $S_n/G$. 
Then we determine those $\alpha\in\A_n$ such that $V(\alpha)$ is $G$-invariant. This yields subgroups $\A_n^G, \TT_n^G$.  
Since the transcendental lattice of a surface is a birational invariant, we obtain
\begin{eqnarray}\label{eq:T^G}
T(X) = \bigoplus_{\alpha\in\TT_n^G} V(\alpha).
\end{eqnarray}
We list the subgroups $\TT_n^G$ as $(\Z/n\Z)^*$-orbits of a single element $\alpha\in\TT_n$. The element is represented by the triple $(a_1, a_2, a_3)$ corresponding to the affine chart (\ref{eq:aff}).
$$
\begin{array}{c|c||c||c}
\hline
m & \lambda & n & \TT_n^G = \mbox{orbit}(\alpha)\\
\hline \hline
4 & 0 & 8 & [4, 2, 1]\\
\hline
8 & 0 & 16 & [8, 5, 1]\\
& \sqrt{3} & 24 & [12, 8, 3]\\
\hline
16 & 0 & 32 & [16, 9, 5]\\
& -2 & 16 & [8, 2, 5]\\
& \sqrt{3} & 48 & [24, 16, 3]\\
\hline
\end{array}
$$
The decomposition (\ref{eq:T^G}) carries over to the Galois representation $\varrho$ associated to $T(X)$. Hence we can compute the zeta function of $X$. Again we refer to the models given by (\ref{eq:D-2}) for $X_{-2}$ and~(\ref{eq:D-3}) for $X_{\sqrt{3}}$.  

\begin{Theorem}\label{Thm:zeta}
Let $X=X_0, X_2$ or  $X_{\sqrt{3}}$ in one of the non-unimodular families except for $X_2, X_{\sqrt{3}}$ in case $m=4$. 
Then
\[
\zeta(X/\F_q, T) = \dfrac 1{(1-T)\, P(T)\, (1-q^2\, T)}
\]
where
\[
P(T) = (1-q\,T)^{22-\phi(n)}  \prod_{\alpha\in\TT_n^G} (1-\jj(\alpha)\,T).
\]
\end{Theorem}

\emph{Proof:}
The product in $P(T)$ is the reciprocal characteristic polynomial of Frobenius on $\varrho$ by (\ref{eq:T^G}). The other factor of $P(T)$ comes from $\NS(X_\C)$. By the above considerations, $\rho(X_\C)=22-\mbox{rank}(T(X))=22-\phi(n)$. Since each fibre component is defined over $\Q$, $\NS(X_\C)$ is generated by divisors over $\Q$. Hence Frobenius operates as multiplication by $q$. \qed

\subsection{Mirror symmetry}
\label{s:mirror}

Mirror symmetry is supposed to interchange complex and K\"ahler structure. For K3 surfaces, we can impose further conditions on the lattices of algebraic and transcendental lattices. Here we employ the
notion of mirror symmetry introduced by Dolgachev \cite{Dol}.

\begin{Definition}
Let $X$ be an algebraic K3 surface. A K3 surface $\breve X$ is a {\em mirror} of $X$ if
\begin{eqnarray}\label{eq:mirror}
T_X = U \oplus S_{\breve X}.
\end{eqnarray}
\end{Definition}

Mirror symmetry is exhibited for families of K3 surfaces. For instance, for the two families of K3 surfaces with $m=4$ in Theorem~\ref{thm}, the mirror family would be general elliptic surfaces with section (so that $\NS=U$) resp.~with bisection (so that $\NS=U(2)$). 

In \cite{LSY}, it is shown that the special member $X_0$ in the unimodular family has mirror surfaces of CM type. This instance of arithmetic mirror symmetry is our motivation to study the familes for $m=8$ and $16$ from Theorem~\ref{thm}. 

Consider the families of K3 surfaces for $m=8, 16$ in Theorem~\ref{thm}. By definition, their general members are mirrors of each other. Here we want to point out that mirror symmetry extends to specific members in an arithmetic way:

\begin{itemize}
\item
The surfaces at $\lambda=0, \pm\sqrt{3}$ have general $\rho$; they are all of CM type.
\item
The surfaces at $\lambda=\pm 2$ degenerate with $T(X)$ of rank $m/2$ instead of $m$. Both have CM type. 
\end{itemize}

In fact, both families of K3 surfaces can be collected in a single family of elliptic surfaces over $\PP^1$. For this we only have to apply the base change 
\[
t\mapsto t^{32/m}.
\]
The resulting elliptic surface $Y$ of Euler number $e(Y)=36$ is given by the following Weierstrass equation:
\[
Y:\;\;\; y^2 = x^3 + \lambda\,x^2 + x + t^{16}.
\]
It has discriminant
\[
\Delta = 16\, (27\,t^{32}-2\,\lambda\,(2\,\lambda^2+9)\,t^{16}-\lambda^2+4),
\]
so in general there are 32 fibres of type $I_1$ plus one fibre of type $IV$ at $\infty$. The degeneration behaviour is the same as before. The general surface in this family has $\rho=10$, since the Mordell-Weil group has rank six. Up to finite index, it is obtained by base change from the family of rational elliptic surfaces
\[
Z:\;\;\; y^2 = x^3 + \lambda\,x^2 + x + t^{4}.
\]
Here the general member has $MW(Z)=E_6^\vee$. It follows that the general member for $Y$ has transcendental lattice of rank 24. Note, however, that despite the non-symplectic automorphism $t\mapsto\zeta_{16}\,t$, the rank of $T(Y)$ is not always divisible by $8$. E.g.~the surface at $\lambda=2$ is of CM type with $\mbox{rank}(T(Y))=12$ by construction.

\vspace{0.5cm}

\textbf{Acknowledgement:} The author would like to thank S.~Kond\= o for stimulating discussions and the referee for many helpful comments. 
Funding from DFG under research grant Schu 2266/2-2 is gratefully acknowledged. 

\vspace{0.5cm}

\vspace{0.8cm}

Matthias Sch\"utt\\
Institute for Algebraic Geometry\\
Leibniz University Hannover\\
Welfengarten 1\\
30167 Hannover\\
Germany\\
schuett@math.uni-hannover.de


\begin{thebibliography}{992}


\bibitem{Dol} I.~Dolgachev, 
\emph{Mirror symmetry for lattice polarized K3 surfaces},
J.~Math.~Sci.~{\bf 81} (1996), 2599--2630. 

%

\bibitem{Inose} H.~Inose,
\newblock Defining equations of singular $K3$ surfaces and a notion of isogeny. 
\newblock In  {\em Proceedings of the International Symposium on Algebraic Geometry} (Kyoto Univ., Kyoto, 1977),  Kinokuniya Book Store, Tokyo (1978), 495--502.

\bibitem{Ko}
S.~Kond\=o, 
\emph{Automorphisms of algebraic K3 surfaces which act trivially
on Picard groups}, 
J.~Math.~Soc.~Japan, {\bf 44}, No.~1 (1992), 75--98.



\bibitem{L} R.~Livn\'e, 
\emph{Motivic orthogonal two-dimensional representations of  Gal$(\overline{\Q}/\Q)$},
Israel J.~of Math.~{\bf 92} (1995), 149--156.

\bibitem{LSY}
R.~Livn\'e, M.~Sch\"utt and N.~Yui,
\emph{The modularity of K3 surfaces with non-symplectic group actions},
preprint (2009), arXiv: 0904.1922.


\bibitem{MO}
N.~Machida and K.~Oguiso,
\emph{K3 surfaces admitting finite non-symplectic group actions},
J.~Math.~Sci.~Univ.~Tokyo, {\bf 5} (1998), No.~2, 273--297. 


\bibitem{N} V.~V.~Nikulin,
\emph{Integral Symmetric Bilinear Forms and Some of their Applications},
Math.~USSR Izvestija, {\bf 14} (1980), No.~1, 103--167.

\bibitem{Ni}
V.~V.~Nikulin, 
\emph{Finite groups of automorphisms of K\"ahler K3 surfaces},
Trans.~Moscow Math.~ Soc.~{\bf 38} (1980), 71--135.)

\bibitem{N-2}
V.~V.~Nikulin,
\emph{On the quotient groups of the automorphism group of hyperbolic forms by the subgroup generated by 2-reflections},
J.~Soviet Math.,~{\bf 22} (1983), 1401--1476.


\bibitem{OZ}
K.~Oguiso and D.-Q.~Zhang,
\emph{On Vorontsov's theorem on K3 surfaces with non-symplectic group-actions},
 Proc.~AMS {\bf 128} (2000), No.~6, 1571--1580.

\bibitem{S-CM}
M.~Sch\"utt,
 \emph{CM newforms with rational coefficients}, 
Ramanujan J.~{\bf 19} (2009), 187--205. 

\bibitem{Sh}
T.~Shioda,
\emph{An explicit algorithm for computing the Picard number of certain algebraic surfaces},
Amer.~J.~Math.~{\bf 108}, No. 2 (1986), 415--432.


\bibitem{Sandwich} T.~Shioda, \emph{Kummer sandwich theorem of certain elliptic K3 surfaces}, Proc.~Japan~Acad.,~{\bf 82}, Ser.~A (2006), 137--140.

\bibitem{Tate} J.~Tate, \emph{Algorithm for determining the type of a singular fibre in an elliptic pencil}, in: \emph{Modular functions of one variable IV} (Antwerpen 1972), SLN {\bf 476} (1975), 33--52.


\bibitem{Vo}
S.~P.~Vorontsov,
\emph{Automorphisms of even lattices arising in connection with automorphisms
of algebraic {$K3$} surfaces} (Russian), 
Vestnik Moskov.~Uni.~Ser.~I mat.~Mekh.~(1983), No.~2, 19--21. 
%

\bibitem{Weil}
A.~Weil,
\emph{Numbers of solutions of equations in finite fields},
Bull.~AMS {\bf 55} (1949), 497--508.



\end{thebibliography}
\end{document}